\documentclass[12pt]{article}
\usepackage{amssymb,amsfonts}%amsLatex

\newcommand{\binom}[2]{{\textstyle \left( \! \! \! \begin{array}{c} {#1} \\ {#2} \end{array} \! \! \! \right)}}
\newcommand{\ona}{\mathrm}
\newcommand{\De}{\Delta}
\newcommand{\Ga}{\Gamma}
\newcommand{\ga}{\gamma}

\newcommand{\va}{\varphi}
\newcommand{\wt}{\widetilde}
\newcommand{\sub}{\subset}
\newcommand{\ti}{\times}
\newcommand{\ot}{\otimes}
\newcommand{\op}{\oplus}
\newcommand{\bw}{\wedge}
\newcommand{\C}{\mathbb{C}}
\newcommand{\Z}{\mathbb{Z}}
\newcommand{\gk}{\mathfrak{g}}
\newcommand{\bno}{\bigskip\noindent}

\newcommand{\Cur}{\mathop{\rm Cur}}
\newcommand{\Mat}{\mathop{\rm Mat}}

\title{Formal Distribution Algebras and Conformal Algebras}

\author{Victor G. Kac \thanks{%
    Department of Mathematics, M.I.T., 
    Cambridge, MA 02139, 
    $<$kac@math.mit.edu$>$
}}

\date{}

\begin{document}

\maketitle

%\begin{abstract}
%\end{abstract}

\section*{Introduction}

Conformal algebra is an axiomatic description of the operator product
expansion (or rather its Fourier transform) of chiral fields in a
conformal field theory. It turned out to be an adequate tool for the
realization of the program of the study of Lie (super)algebras and
associative algebras (and their representations), satisfying the sole
locality property \cite{K3}.
The first basic definitions and results appeared in my book
\cite{Kac}  and review \cite{K3}.
In the present paper I review recent developments in conformal
algebra, including some of \cite{Kac} and \cite{K3} but in a different
language. Here I use the $\lambda$-product, which is the Fourier transform
of the OPE, or, equivalently, the generating series of the $n$-th
products used in  \cite{Kac} and \cite{K3}. This makes the exposition
much more elegant and transparent.

Most of the work has been done jointly with my collaborators. The
structure theory of finite Lie conformal algebras is a joint
paper with A. D'Andrea \cite{DK}. The theory of conformal modules
has been developed with S.-J. Cheng \cite{CK} and of their
extensions with S.-J. Cheng and M. Wakimoto \cite{CKW}. The
understanding of conformal algebras $ \ona{Cend}_N $ and $gc_N$
was achieved with A. D'Andrea \cite{DK}, and of their finite
representations with B.~Bakalov, A.~Radul and M.~Wakimoto
\cite{BKRW}. The connection to $\Ga$-local and $\Ga$-twisted
formal distribution algebras has been established with M.
Golenishcheva-Kutuzova \cite{GK} and with B.~Bakalov and A.~D'Andrea
\cite{BDK}.  Cohomology theory has been worked out with B.~Bakalov and A.~Voronov \cite{BKV}.

\section{Calculus of formal distributions}

Let $U$ be a vector space over $\C.$ A $U$-valued \emph{formal
distribution} in one indeterminate $z$ is a linear $U$-valued function
on the space of Laurent polynomials $\C[z, z^{-1}].$ Such a formal
distribution $a(z)$ can be uniquely written in the form
$a(z)=\sum_{n\in\Z}a_nz^{-n-1},$
where $a_n\in U$ is defined by $a_n=\ona{Res}_zz^na(z)$ and
$\ona{Res}_z$ stands for the coefficient of $z^{-1}$. The space of
these distributions is denoted by $U[[z, z^{-1}]].$

Likewise, one defines a $U$-valued formal distribution in $z$ and $w$
as a linear function on the space of Laurent polynomials in $z$ and
$w,$ and such a formal distribution can be uniquely written in the
form
$a(z, w)=\sum_{m, n\in\Z}a_{m, n}z^{-m-1}w^{-n-1}.$
This formal distribution defines a linear map
$
  D_{a(z, w)}: {\C}[w, w^{-1}]\!\rightarrow U[[w, w^{-1}]]
$
by letting
$(D_{a(z, w)}f)(w)\!=\!\ona{Res}_za(z, w)f(z).$
The most important $\C$-valued formal distribution in $z$ and $w$ is
the \emph{formal} $\delta$-\emph{function} $\delta(z-w)$ defined by
$D_{\delta(z-w)}f=f.$
Explicitly:
$$
  \delta(z-w)=z^{-1}\sum_{n\in\Z}(w/z)^n.
$$

A formal distribution $a(z, w)$ is called \emph{local} if
$(z-w)^Na(z, w)=0$ for $N\gg 0.$
Note that $a(w, z),$ $\partial_za(z, w)$ and $\partial_wa(z, w)$ are local if
$a(z, w)$ is. It is easy to show (\cite{Kac}, Corollary 2.2) that a
formal distribution $a(z, w)$ is local iff it can be
represented as a finite sum of the form
\begin{equation}
a(z, w)=\sum_{j\in\Z_+}c^j(w)\partial_w^{(j)}\delta(z-w).
\end{equation}
(Such a representation is unique.)
Here and further $\partial^{(j)}$ stands for $\partial^j/j!.$
This is called the \emph{operator product expansion} (OPE), the
$c^j(w),$ called the OPE coefficients, being given by
\begin{equation}
c^j(w)=\ona{Res}_za(z, w)(z-w)^j.
\end{equation}
Note that $a(z, w)$ is local iff $D_{a(z, w)}$ is a differential
operator:
\hfill \\
$D_{a(z, w)}=\sum_{j\in\Z_+}c^j(w)\partial_w^{(j)}.$
Note also that $D_{a(w, z)}$ is the adjoint differential operator:
$D_{a(w, z)} =\sum_{j\in\Z_+}(-\partial_w)^{(j)}c^j(w).$

In order to study the properties of the expansion (1), it is
convenient to introduce the \emph{formal Fourier transform} of a formal
distribution $a(z, w)$ by the formula:
$$\Phi_{z, w}^{\lambda}(a(z, w))=\ona{Res}_ze^{\lambda(z-w)}a(z, w).$$
This is a linear map from $U[[z, z^{-1}, w, w^{-1}]]$ to $U[[\lambda]]
[[w, w^{-1}]].$ We have:
$$\Phi^{\lambda}_{z, w}(\partial^j_w\delta(z-w)) =\lambda^j.$$
Hence the formal Fourier transform of the expansion (1) is
\begin{equation}
\Phi_{z, w}^{\lambda}(a(z, w))= \sum_{n\in\Z_+}\lambda^{(n)}c^n(w).
\end{equation}
(As before, $\lambda^{(n)}$ stands for $\lambda^n/n!.$) In other words, the
formal Fourier transform of a formal distribution $a(z, w)$ is the
generating series of its OPE coefficients. If $a(z, w)$ is
local then its formal Fourier transform is polynomial in $\lambda.$

We have the following important relations:
\begin{eqnarray}
  \Phi^{\lambda}_{z, w}\partial_z
  & = &
  -\lambda\Phi^{\lambda}_{z, w}
  =
  \left[\partial_w, \Phi^{\lambda}_{z, w}\right],
  \\
  \Phi^{\lambda}_{z, w}a(w, z)
  & = &
  \Phi^{-\lambda-\partial_w}_{z, w}a(z, w)
  \quad
  \mathrm{if~} a(z, w) \mathrm{~is~local}.
\end{eqnarray}
(The right-hand side of (5) means that the indeterminate $\lambda$ in (3)
is replaced by the operator $-\lambda-\partial_w.$)

\bigskip\noindent
\textbf{Remark 1.}\quad Formulas (4) and (5) are equivalent to the
following relations for the OPE coefficients $c^n_z(w),$ $c^n_w(w)$ and
$\wt{c}^n(w)$ of the formal distributions $\partial_za(z, w)$,
$\partial_wa(z, w)$ and $a(w, z)$ respectively:
$$
  \begin{array}{c}
    c^n_z(w) = -nc^{n-1}(w),
    \quad
    c^n_w(w) = \partial_wc^n(w)+ nc^{n-1}(w),
    \\
    \quad
    \wt{c}^n(w) = \sum_{j\in\Z_+}(-1)^{j+n}\partial^{(j)}_wc^{n+j}(w).
  \end{array}
$$

A composition of two Fourier transforms, like
$\Phi^{\lambda}_{z, w}\Phi^{\mu}_{x, w}$, is a linear map from
%\linebreak
\hfill \\
$U[[z, z^{-1}, w, w^{-1}, x, x^{-1}]]$ to
$U[[\lambda, \mu]][[w, w^{-1}]].$ The following relation is
of fundamental importance:
\begin{equation}
\Phi^{\lambda}_{z, w}\Phi^{\mu}_{x, w}
= \Phi^{\lambda+\mu}_{x, w}\Phi^{\lambda}_{z, x}\, .
\end{equation}
The proof is very easy:
$$
  \ona{Res}_z\ona{Res}_x e^{\lambda(z-w)+\mu(x-w)}a(z, w, x)
  =
  \ona{Res}_z\ona{Res}_x e^{\lambda(z-x)}e^{(\lambda+\mu)(x-w)} a(z, w, x).
$$

\section{Formal distribution algebras}
\label{sec:2}
Now let $U$ be an algebra. Given two $U$-valued formal distributions
$a(z)$ and $b(z),$ we may consider the formal distribution
$
  a(z)b(w) = \sum_{m, n}a_mb_nz^{-m-1}w^{-n-1}
$ and its formal Fourier
transform, which we denote by
$$a(w)_{\lambda}b(w)\equiv
\sum_{j\in\Z_+}\lambda^{(j)}(a(w)_{(j)}b(w))
=\Phi^{\lambda}_{z, w}(a(z)b(w)).$$
This is called the $\lambda$-\emph{product}. The coefficients $a(w)_{(j)}b(w),$
$j\in\Z_+,$ called $j$-th products, can be calculated by (2).

A pair of $U$-valued
formal distributions $a(z)$ and $b(z)$ is called
\emph{local} if the formal distribution $a(z)b(w)$ is local (this is
not symmetric in $a$ and $b$ in general). A subset $F\sub U[[z,
z^{-1}]]$ is called a \emph{local family} if all pairs of
formal distributions
from $F$ are local. We denote by $\bar{F}$ the minimal
subspace of $U[[z, z^{-1}]]$ which contains $F$ and is closed under
all $j$-th products.

\bno
\textbf{Lemma 2 } \cite{Kac}\textbf{.} \emph{
Let $U$ be a Lie (super)algebra or an associative algebra, and let $F\sub
U[[z, z^{-1}]]$ be a local family. Then $\bar{F}$ is a local family as
well.
}

\bigskip
A pair $(U, F),$ where $U$ is an algebra and $F\sub U[[z, z^{-1}]],$ is
called a \emph{formal distribution algebra} if $\bar{F}$ is a local
family of formal distributions whose coefficients span $U$ over $\C.$

Let $\ona{Con}(U, F) = \C[\partial_z]\bar{F}.$
This is still a local family, hence the $\lambda$-product defines a map:
$
  \ona{Con}(U, F)\ot_{\C}\ona{Con}(U, F)\rightarrow
  \C[\lambda]\ot_{\C}\ona{Con}(U, F)
$.
Due to (4), the $\lambda$-product satisfies
\begin{equation}
  \begin{array}{rcl}
    \partial_za(z)_{\lambda}b(z)
    & = &
    -\lambda a(z)_{\lambda}b(z),
%  \tag{7a}
    \\
    a(z)_{\lambda}\partial_zb(z)
    & = &
    (\partial_z+\lambda)(a(z)_{\lambda}b(z)).
  \end{array}
%  \tag{7b}
  \label{Kac7}
\end{equation}

\bno
\textbf{Example 2.}\quad Let $A$ be an algebra and
% let
$U=A[t, t^{-1}]$,
\hfill \\
$ F = \{a(z)=\sum_{n\in\Z}(at^n)z^{-n-1}\}_{ a\in A}$.
Since $a(z)b(w)= (ab)(w)\delta(z-w)$ for any $a, b\in A,$ we see that
$F=\bar{F}$ is a local family, hence $(U, F)$ is a formal distribution
algebra. It is called a \emph{current algebra}.

\bno
\textbf{Remark 2.1.}\quad For an algebra $U$ denote by $U^{op}$ the
algebra with the opposite multiplication $\circ.$ It follows from (5)
that if $(U, F)$ is a formal distribution algebra with $\lambda$-product
$a(z)_{\lambda}b(z),$ then $(U^{op}, F)$ is a formal distribution
algebra with $\lambda$-product
\begin{equation}
  a(z)_{\lambda}\circ b(z) =b(z)_{-\lambda-\partial_z}a(z).
\end{equation}

In the case when $U$ is a Lie algebra we shall denote the
$\lambda$-product by $[a(z)_{\lambda}b(z)]$ and call it the $\lambda$-{\it
bracket}.

\bno
\textbf{Proposition 2.}\quad
{\it
Let $(U, F)$ be a formal distribution algebra.
\begin{itemize}
\item[{\em (a)}]
$U$ is an associative algebra iff
\begin{equation}
a(z)_{\lambda}\left(b(z)_{\mu}c(z)\right)
= \left(a(z)_{\lambda}b(z)\right)_{\lambda+\mu}c(z).
\end{equation}
\item[{\em (b)}]
$U$ is a commutative algebra iff
\begin{equation}
a(z)_{\lambda}b(z)= b(z)_{-\lambda-\partial_z}a(z).
\end{equation}
\item[{\em (c)}]
$U$ is a Lie (super)algebra iff
\begin{equation}
  \begin{array}{rcl}
    \left[a(z)_{\lambda}b(z)\right]
    & = &
    -p(a, b)\left[b(z)_{-\lambda-\partial_z}a(z)\right],
%  \tag{11a}
    \\
    \left[a(z)_{\lambda}\left[b(z)_{\mu}c(z)\right]\right]
    & = &
    \left[ \left[ a(z)_{\lambda} b(z)\right]_{\lambda + \mu}  c(z) \right]
    +
    p(a, b)\left[b(z)_{\mu}\left[a(z)_{\lambda}c(z)\right]\right]
  \end{array}
  \label{Kac11}
%  \tag{11b}
\end{equation}
($p(a, b)$ stands for $(-1)^{p(a)p(b)},$ where $p(a)$ is the parity of
an element $a\in U$.)
\end{itemize}
}

\bno
\emph{Proof.}\quad It follows from (5) and (6).

\bno
\textbf{Remark 2.2.}\quad Let $(U, F)$ be a $\Z_2$-graded formal
distributions associative algebra where the family $F$ is compatible
with the $\Z_2$-gradation. Then $(U_-, F),$ where $U_-$ is the Lie
superalgebra associated to $U$ with the bracket $[a, b]=ab-p(a, b)ba,$
is a formal distribution Lie superalgebra with the $\lambda$-bracket
$$[a(z)_{\lambda}b(z)]
=a(z)_{\lambda}b(z)-p(a, b)b(z)_{-\lambda-\partial}a(z).$$

\bno
\textbf{Remark 2.3.}\quad (Cf. Remark 1.)
$\ona{Con}(U, F)$ is a $\C[\partial_z]$-module with a $\C$-bilinear product
$a(z)_{(j)}b(z)$ for each $j\in\Z_+$ such that:
\begin{eqnarray*}
  a(z)_{(j)}b(z)
  & = &
  0\qquad\mathrm{for~} j\gg 0,\\
  \partial_za(z)_{(j)}b(z)
  & = &
  -ja(z)_{(j-1)}b(z),\\
  a(z)_{(j)}\partial_zb(z)
  & = &
  \partial_z(a(z)_{(j)}b(z))+ja(z)_{(j-1)}b(z).
\end{eqnarray*}
Algebra $U$ is commutative iff
$$a(z)_{(n)}b(z)=\sum_{j\in\Z_+}(-1)^{j+n}\partial_z^{(j)}
(b(z)_{(n+j)}a(z)).$$
Algebra $U$ is associative iff
$$
  a(z)_{(m)}(b(z)_{(n)}c(z))
  =
  \sum_{j\in\Z_+}
    \binom{m}{j}
    \left(a(z)_{(j)}b(z)\right)_{(m+n-j)}c(z).
$$
Algebra $U$ is a Lie superalgebra iff
\begin{eqnarray*}
  a(z)_{(n)}b(z)
  & = &
 - p(a, b)
  \sum_{j\in\Z_+}(-1)^{j+n}\partial_z^{(j)}
  (b(z)_{(n+j)}a(z)),
  \\
  \left[a(z)_{(m)}, b(z)_{(n)}\right] c(z)
  & = &
  \sum_{j\in\Z_+}
    \binom{m}{j}
    (a(z)_{(j)}b(z))_{(m+n-j)}c(z).
\end{eqnarray*}

\section{Conformal algebras}

Motivated by the discussion in the previous Section, we give the
following definitions \cite{Kac}:

A \emph{conformal algebra} is a $\C[\partial]$-module $R$ endowed with the
$\lambda$-product $a_{\lambda}b$ which defines a linear map $R\ot_{\C}R\rightarrow \C[\lambda]
\ot_{\C}R$ subject to the following axiom (cf. (\ref{Kac7})):
\begin{equation}
(\partial a)_{\lambda}b=-\lambda a_{\lambda}b,\quad
a_{\lambda}\partial b=(\partial+\lambda)(a_{\lambda}b).
\end{equation}
We write $a_{\lambda}b=\sum_{j\in\Z_+}\lambda^{(j)}(a_{(j)}b)$ and call
$a_{(j)}b$ the $j$-\emph{th product} of $a$ and $b.$

A conformal algebra is called \emph{associative} if (cf. (9)):
\begin{equation}
a_{\lambda}(b_{\mu}c)= (a_{\lambda}b)_{\lambda+\mu}c,
\end{equation}
and \emph{commutative} if (cf. (10)):
\begin{equation}
a_{\lambda}b=b_{-\lambda-\partial}a.
\end{equation}
The skewcommutativity and Jacobi identity read (cf. (\ref{Kac11})):
\begin{equation}
  \begin{array}{rcl}
    \left[a_{\lambda}b\right]
    & = &
    -p(a, b)\left[b_{\lambda}a\right],
%  \tag{15a}
    \\
    \left[a_{\lambda}\left[b_{\mu}c\right]\right]
    & = &
    \left[[a_{\lambda}b]_{\lambda+\mu}c\right]+p(a, b)[b_{\mu}[a_{\lambda}c]].
  \end{array}
%  \tag{15b}
  \label{Kac15}
\end{equation}
A $\Z_2$-graded conformal algebra with the $\lambda$-bracket satisfying
(\ref{Kac15}) is called a \emph{Lie conformal superalgebra}. A $\Z_2$-graded
associative conformal algebra with the $\lambda$-bracket (cf. Remark 2.2):
\begin{equation}
[a_{\lambda}b]=a_{\lambda}b-p(a, b)b_{-\lambda-\partial}a
\end{equation}
is a Lie conformal superalgebra.

All these notions and formulas can be given in terms of $n$-th
products (cf. Remark 2.3).

\bno
\textbf{Remark 3.}\quad If $R$ is a conformal algebra, then $\partial R$ is
its ideal with respect to the 0-th product $a_{(0)}b=
a_{\lambda}b|_{\lambda=0}.$

\bno
\textbf{Proposition 3} \cite{DK}\textbf{.}\quad {\it
Any torsion element of a conformal algebra is central.
}

\bno
\emph{Proof}.\quad If $P(\partial)a=0$ for some $P(\partial)\in\C[\partial],$ we have:
$$
  0=(P(\partial)a)_{\lambda}b=P(-\lambda)a_{\lambda}b, \;
  0=b_{\lambda}P(\partial)a=P(\partial+\lambda)
(b_{\lambda}a).
$$
It follows that $a_{\lambda}b=0=b_{\lambda}a.$ \hfill  $\square$

\bigskip
As we have seen in the previous section, given a formal distribution
algebra $(U, F),$ one canonically associates to it a conformal algebra
$\ona{Con}(U, F)$. An ideal $I$ of $U$ is called \emph{irregular} if
there exists no non-zero $b(z)\in\bar{F}$ such that all $b_n\in I.$
Denote the image of $F$ in $U/I$ by $F_1.$ It is clear that the
canonical homomorphism $U\rightarrow U/I$ induces a surjective homomorphism
$\ona{Con}(U, F)\rightarrow \ona{Con}(U/I, F_1),$ which is an isomorphism iff the
ideal $I$ is irregular.

Note that $\ona{Con}$ is a functor from the category of formal distribution
algebras with morphisms being homomorphisms $\varphi:(U, F)\rightarrow (U_1, F_1)$
such that $\varphi(F)\sub\bar{F}_1,$ to the category of conformal
algebras with morphisms being all homomorphisms.

In order to construct the (more or less) reverse functor, we need the
notion of \emph{affinization} $\wt{R}$ of a conformal algebra $R$
(which is a generalization of that for vertex algebras \cite{Borcherds}). We
let $\wt{R}=R[t, t^{-1}]$ with $\wt{\partial}=\partial+\partial_t$ and the
$\lambda$-product  \cite{Kac}:
$af(t)_{\lambda}bg(t)=(a_{\lambda+\partial_t}b)f(t)g(t')|_{t'=t}.$
In terms of $k$-th products this formula reads:
\hfill \\
$
  at^m_{(k)}bt^n
  =
  \sum_{j\in\Z_+}
    \binom{m}{j}
    \left(a_{(j+k)}b\right)t^{m+n-j}.$
But, by Remark 3, $\wt{\partial}\wt{R}$ is an ideal of $\wt{R}$ with
respect to the 0-th product. We let $\ona{Alg}R=\wt{R}/\wt{\partial}\wt{R}$
with the 0-th product and let $\Re=
\{\sum_{n\in\Z}(at^n)z^{-n-1}\}_{a\in R}.$ Then $(\ona{Alg}R, \Re)$ is
a formal distribution algebra. Note that Alg is a functor from the
category of conformal algebras to the category of formal distribution
algebras. One has \cite{Kac}:
$$
 \ona{Con}(\ona{Alg}R)
 =
 R,
 \qquad
 \ona{Alg}(\ona{Con}(U, F))
 =
 (\ona{Alg}\bar{F}, \bar{F}).
$$
Note also that $(\ona{Alg}R, \Re)$ is the maximal formal distribution
algebra associated to the conformal algebra $R$, in the sense that all
formal distribution algebras $(U, F)$ with $\ona{Con}(U, F)=R$ are
quotients of $(\ona{Alg}R, \Re)$ by irregular ideals. Such formal
distribution algebras are called \emph{equivalent}.

We thus have an equivalence of categories of conformal algebras and
equivalence classes of formal distribution algebras. This
equivalence restricts to the categories of associative, commutative
and Lie (super)algebras. So the study of formal distribution algebras
reduces to the study of conformal algebras.

\bno
\textbf{Example 3.1.}\quad Let $A$ be an algebra and let $A[t, t^{-1}]$
be the associated current formal distribution algebra (see Example
2). Then the associated conformal algebra is
$\ona{Cur}A=\C[\partial]\ot_{\C}A$
with multiplication defined by
$a_{\lambda}b=ab$ for  $a, b\in A$
(and extended to $\ona{Cur}A$ by (12)). This is called the {\it
current conformal algebra} associated to $A.$ Note that $A[t, t^{-1}]$
is the maximal formal distribution algebra associated to
$\ona{Cur}A,$ and that for any non-invertible Laurent polynomial
$P(t)$, the formal distribution algebra $A[t, t^{-1}]/(P(t))$ is
equivalent to $A[t, t^{-1}].$

\bno
\textbf{Example 3.2.}\quad The simplest formal distribution Lie
algebra beyond current algebras is the Lie algebra
$\ona{Vect}\C^{\ti}$ of regular vector fields on $\C^{\ti}$ ($=$ Lie
algebra of infinitesimal conformal transformation of $\C^{\ti}$, hence
the choice of the term ``conformal''). Vector fields
$L_n=-t^{n+1}\partial_t\; (n\in\Z)$ form a basis of Vect\,$\C^{\ti}$ with
the familiar commutation relation
$[L_m, L_n]=(m-n)L_{m+n}.$ Furthermore,
$L(z)= - \sum_{n\in\Z}(t^n\partial_t)z^{-n-1}$ is a local formal distribution
(i.e. the pair $(L, L)$ is local) since
$$
  [L(z), L(w)]=\partial_wL(w)\delta(z-w)+2L(w)\delta'_w(z-w).
$$
The associated Lie conformal algebra is $\C[\partial]L$ with
$\lambda$-bracket:  $[L_{\lambda}L]=(\partial+2\lambda)L.$ This is called the
\emph{Virasoro conformal algebra} since the centerless Virasoro algebra
Vect$\,\C^{\ti}$ is the maximal (and only) associated formal
distribution Lie algebra.

A formal distribution algebra $(U, F)$ (resp. a conformal algebra $R$)
is called \emph{finite} if the $\C[\partial]$-module
$\C[\partial]\bar{F}$ (resp. $R$) is finitely generated.

\section{Structure theory of finite conformal (super)algebras}

A conformal algebra is called simple if it is not commutative and
contains no nontrivial ideals.

\bno
\textbf{Theorem 4.1} \cite{DK}\textbf{.}\quad {\it
A simple finite Lie conformal algebra is isomorphic either to a
current conformal algebra $\ona{Cur}\gk$, where $\gk$ is a simple
finite-dimensional Lie algebra, or to the Virasoro conformal algebra.
}

\medskip
Of course, translating this into the language of formal distribution
Lie algebras, we obtain the following result:

\bno
\textbf{Corollary 4.}\quad {\it
Any formal distribution Lie algebra which is finite and simple
(i.e. any its non-trivial ideal is irregular) is isomorphic either to
$(\ona{Vect} \C^{\ti}, \{L(z)\})$ or to a quotient of $(\gk[t, t^{-1}],
\{a(z)\}_{a\in\gk})$ where $\gk$ is a simple finite-dimensional Lie
algebra.
}

\bno
\textbf{Open Problem.}\quad {\it
Classify simple formal distribution Lie algebras $(\gk, F)$ for which $F$ is
a finite set.}

\medskip
Under the assumption that $\gk$ is $\Z$-graded the only
possibilities are the Virasoro and twisted loop algebras.  This follows
from a very difficult theorem of Mathieu \cite{M}, but there are
many examples which are not $\Z$-graded.

\medskip

The $\C$-span of all elements of the form $a_{(m)}b$ of a conformal
algebra $R,$ $m\in\Z_+,$ is called the \emph{derived algebra} of $R$
and is denoted by $R'.$ It is easy to see that $R'$ is an ideal of $R$
such that $R/R'$ has a trivial $\lambda$-product. We have the derived
series $R\supset R'\supset R''\supset\cdots.$ A conformal algebra is called
\emph{solvable} if the $n$-th member of this series is zero for $n\gg
0.$ A Lie conformal algebra is called \emph{semisimple} if it contains
no non-zero solvable ideals.

\bno
\textbf{Theorem 4.2} \cite{DK}\textbf{.}\quad {\it
Any finite semisimple Lie conformal algebra is a direct sum of
conformal algebras of the following three types:
\begin{itemize}
\item[{\em (i)}]
current conformal algebra $\ona{Cur}\gk,$ where $\gk$ is a semisimple
finite-dimensional Lie algebra,
\item[{\em (ii)}]
Virasoro conformal algebra,
\item[{\em (iii)}]
the semidirect product of {\em (i)} and {\em (ii)}, defined by
$L_{\lambda}a=(\partial+\lambda)a$ for $a\in\gk.$
\end{itemize}
}

\medskip
The proof of these results uses heavily Cartan's theory of filtered
Lie algebras.

By far the hardest is the classification of finite simple Lie
conformal superalgebras. The list is much richer than that of finite
simple Lie conformal algebras. First, there are many more simple
finite-dimensional Lie superalgebras (classified in \cite{K1}), and
the associated current conformal superalgebras are finite and simple.
Second, there are many ``superizations'' of the Virasoro conformal
algebra. They are associated to superconformal algebras constructed in
\cite{KL} and \cite{CK1}. We describe them below.

Let $\bw(N)$ denote the Grassmann algebra in $N$ indeterminates
$\xi_1, \ldots, \xi_N$ and let $W(N)$ be the Lie superalgebra of all
derivations of the superalgebra  $\bw(N)$. It consists of all linear
differential operators $\sum_iP_i\partial_i,$ where $\partial_i$ stands for the
partial derivative by $\xi_i.$

The first series of examples is the series of Lie conformal
superalgebra $W_N$ of rank $(N+1)2^N$ over $\C[\partial]:$
$W_N=\C[\partial]\ot_{\C}(W(N)+\bw(N))$
with the following $\lambda$-brackets $(a, b\in W(N);\; f, g\in\bw(N)):$
$$
  [a_{\lambda}b] =[a, b],\quad [a_{\lambda}f]=a(f)-p(a, f)\lambda fa,\quad
  [f_{\lambda}g] =-(\partial+2\lambda)fg.
$$
The second series is $S_N$ of rank $N2^N$ over $\C[\partial]:$
$S_N=\{D\in W_N| \ona{div}D=0\},$
where
$$\ona{div}(\sum_iP_i(\partial, \xi)\partial_i +f(\partial, \xi))
=\sum_i(-1)^{p(P_i)}\partial_iP_i-\partial f.$$
The third series is a deformation of $S_N$ ($N$
even):  $\tilde{S}_N=\{ D \in W_N | \ona{div}(1+ \xi_1 \ldots \xi_N)
D=0 \}$.  
The fourth series is $K_N$, which is also a subalgebra of $W_N$ (of
rank $2^N$), but it is more convenient to describe it as follows:
$K_N=\C[\partial]\ot_{\C}\bw(N)$
with the following $\lambda$-bracket $(f, g\in\bw(N)):$
$$
  [f_{\lambda}g]
  =
  \left(\frac{1}{2}|f|-1\right)\partial(fg)
  +\frac{1}{2}(-1)^{|f|}\sum_i(\partial_if)(\partial_ig)\\
  +\lambda\left(\frac{1}{2}|f|+\frac{1}{2}|g|-2\right)fg.
$$
We assume here that $f$ and $g$ are homogeneous elements of degree
$|f|$ and $|g|$ in the $\Z$-gradation $\deg\xi_i=1$ for all $i.$

\bno
\textbf{Theorem 4.3} \cite{K2}, \cite{K4}\textbf{.}\quad {\it
Any simple finite Lie conformal superalgebra is isomorphic either to a
current conformal superalgebra $\ona{Cur}\gk,$ where $\gk$ is a simple
finite-dimensional Lie superalgebra, or to one of the conformal
superalgebras $W_N,$ $S_{N+2},$ $\tilde{S}_{2N+2},$ $K_N\,(N\in\Z_+),$ 
or to the
exceptional conformal superalgebra $CK_6$ of rank 32 (which is a
subalgebra of $K_6$) constructed in \cite{CK1}.
}

\medskip
The structure theory of finite associative conformal algebras is much
simpler than that of Lie conformal algebras. Define the central
series $R\supset R^1\supset R^2\supset\cdots$ by letting $R^1=R'$ and
$R^n=\C$-span of all products $a_{(j)}b$ and $b_{(j)}a,$ $j\in\Z_+,$
where $a\in R,$ $b\in R^{n-1},$ for $n\geq 2.$ A conformal algebra is
called \emph{nilpotent} if $R^n=0$ for $n\gg 0.$ An associative
conformal algebra is called semisimple if it contains no non-zero
nilpotent ideals.

\bno
\textbf{Theorem 4.4.}\quad \emph{
Any finite semisimple associative conformal algebra is a direct sum of
associative conformal algebras of the form $\ona{Cur}(\ona{Mat}_N),$
$N\geq 1,$ where $\ona{Mat}_N$ stands for the associative algebra of
all complex $N\ti N$ matrices.
}

\section{Conformal modules and modules over conformal algebras}

Now let $A$ be an associative algebra or a Lie (super)algebra, and let
$V$ be an $A$-module. Given an $A$-valued formal distribution $a(z)$
and a $V$-valued formal distribution $v(z)$ we may consider the formal
distribution $a(z)v(w)$ and its formal Fourier transform
$a(w)_{\lambda}v(w)=\Phi^{\lambda}_{z, w}(a(z)v(w)).$
This is called the $\lambda$-\emph{action} of $A$ on $V.$ It has properties
similar to (\ref{Kac7}):
\begin{equation}
  \begin{array}{rcl}
    \partial_za(z)_{\lambda}v(z)
    & = &
    -\lambda a(z)_{\lambda}v(z),
%  \tag{17a}
    \\
    a(z)_{\lambda}\partial_zv(z)
    & = &
    (\partial_z+\lambda)(a(z)_{\lambda}b(z)).
  \end{array}
%  \tag{17b}
  \label{Kac17}
\end{equation}
In the case when $A$ is associative we have a formula similar to (9):
\begin{equation}
a(z)_{\lambda}(b(z)_{\mu}v(z))
=(a(z)_{\lambda}b(z))_{\lambda+\mu}v(z).
\end{equation}
Likewise, in the case when $A$ is a Lie (super)algebra, we have a
formula similar to (\ref{Kac11}b):
\begin{equation}
[a(z)_{\lambda}, b(z)_{\mu}]v(z)
=[a(z)_{\lambda}b(z)]_{\lambda+\mu}v(z).
\end{equation}
As before, the pair $(a(z),\,v(z))$ is called \emph{local} if the formal
distribution $a(z)v(w)$ is local.

\bno
\textbf{Lemma 5} \cite{K3}\textbf{.}\quad {\it
Let $F \sub A[[z, z^{-1}]]$ be a local family and let $E\sub V[[z,
z^{-1}]]$ be such that all pairs $(a(z), v(z)),$ where
$a(z)\in F$ and $v(z)\in E,$ are local.
Let $\bar{E}$ be the minimal subspace of $V[[z,
z^{-1}]]$ which contains $E$ and all $a(z)_{(j)}v(z)$ for $a(z)\in
F,$ $v(z)\in \bar{E}.$ Then all pairs $(a(z), v(z))$ with
$a(z)\in \bar{F}$ and
$v(z)\in\bar{E}$ are local. Moreover
$a(z)_{(j)}(\C[\partial_z]\bar{E})\sub\C[\partial_z]\bar{E}$ for all
$a(z)\in\C[\partial_z]\bar{F}.$
}

\medskip
Let $F$ be a local family that spans $A$ and let $E\sub
V[[z, z^{-1}]]$ be a family that spans $V.$ Then $(V, E)$ is called
a \emph{conformal module} over the formal distribution algebra $(A,
F)$ if all pairs $(a(z), v(z)),$ where
$a(z)\in F$ and $v(z)\in E$ are local. It follows
from Lemma 5 that a conformal module $(V, E)$ over a formal
distribution (associative or Lie) algebra gives rise to a {\it
module} $\ona{Con}(V, E):=\C[\partial_z]\bar{E}$ over the conformal algebra
$\ona{Con}(A, F).$

A conformal module $(V, E)$ is called \emph{finite} if $\ona{Con}(V, E)$
is a finitely generated $\C[\partial_z]$-module and is called {\it
irreducible} if it contains no irregular submodules.

The definition of a module over a conformal algebra (associative or
Lie) is motivated by (\ref{Kac17}-19) and is given along the same lines as
before. A \emph{module} over a conformal algebra $R$ is a $\C[\partial]$-module
$M$ endowed with the $\lambda$-action $a_{\lambda}v$ which defines a map
$R\ot_{\C} M\rightarrow\C[\lambda]\ot_{\C}M$ subject to the following axioms:
\begin{eqnarray}
\label{eq:20}
  (\partial a)_{\lambda} v
  & = &
  -\lambda a_{\lambda}v,
  \quad
  a_{\lambda}\partial v=(\partial+\lambda)(a_{\lambda} v),  \\
\label{eq:21}
  a_{\lambda}(b_{\mu}v)
  & = &
  (a_{\lambda}b)_{\lambda+\mu}v
  \quad
  \mathrm{if~} R \mathrm{~is~associative}, \\
\label{eq:22}
  \left[a_{\lambda}, b_{\mu}\right]v
  & = &
  \left[a_{\lambda}b\right]_{\lambda+\mu}v
  \quad
  \mathrm{if~} R \mathrm{~is~Lie.}
\end{eqnarray}
A 1-dimensional (over $\C$) $\C[\partial]$-module over a conformal algebra
$R$ is called \emph{trivial} if $a_{\lambda}=0$ for all $a\in R.$

\bno
\textbf{Remark 5.1.}\quad It is easy to see that
$(b_{-\lambda-\partial}a)_{\mu}v=(b_{\mu-\lambda}a)_{\mu}v.$
It follows that a module over a ($\Z_2$-graded) associative conformal
algebra is a module over the corresponding Lie conformal
(super)algebra defined by (16).

In the same way as in Section 3 we have an equivalence of categories
of modules over an associative (resp. Lie) conformal algebra $R$ and
equivalence classes of conformal modules over the associative
(resp. Lie) algebra $\ona{Alg}R.$ Proof of the following proposition
is the same as that of Proposition 3.

\bno
\textbf{Proposition 5} \cite{DK}\textbf{.}\quad {\it
Let $M$ be a module over a conformal algebra $R.$ Then the torsion of
$R$ acts trivially on $M$ and $R$ acts trivially on the torsion of
$M.$
}

\bno
\textbf{Example 5.1.}\quad Let $A$ be an associative or Lie (super)algebra
and let $U$ be an $A$-module. Then, in the obvious way, $(U[t,
t^{-1}], E)$ is a conformal module over the current algebra (see
Example 2.1) with $E =\{u(z)=\sum_n(ut^n)z^{-n-1}\}_{u\in U}.$ The
associated module over the current conformal algebra $\ona{Cur}A$ is
$M(U)=\C[\partial]\ot_{\C}U$
with the $\lambda$-action
$a_{\lambda}u=a(u),\quad a\in A,\; u\in U.$

\bno
\textbf{Example 5.2.}\quad For each $\De,$ $\alpha\in\C$
define the space of densities:
\hfill \\
$F(\De, \alpha)=\C[t, t^{-1}]e^{-\alpha t}(dt)^{1-\De}.$
This is a conformal module over $\ona{Vect}\C^{\ti}$ with
$ E =\{m(z)=\sum_n(t^ne^{-\alpha t}(dt)^{1-\De})z^{-n-1}\}.$ It is
irreducible unless $\De=0$ (in this case it has a regular submodule
$d(\C[t, t^{-1}]e^{-\alpha t})$). The associated module over the Virasoro
conformal algebra (see Example 3.2) is
$M(\De, \alpha)=\C[\partial]m$
with the $\lambda$-action
$L_{\lambda}m=(\partial+\alpha+\De\lambda)m.$

\bno
\textbf{Theorem 5.1} \cite{CK}\textbf{.}\quad

{\it
{\em (a)}
Any non-trivial irreducible finite module over the Virasoro conformal
algebra is isomorphic to $M(\De, \alpha)$ with $\De\neq 0.$

{\em (b)}
Any non-trivial irreducible finite module over the current conformal algebra
$\ona{Cur}\gk,$ where $\gk$ is a finite-dimensional semisimple Lie
algebra, is isomorphic to $M(U)$ where $U$ is a non-trivial
irreducible finite-dimensional $\gk$-module.
}

\bno
\textbf{Corollary 5.}\quad

{\it
{\em (a)}
Any non-trivial irreducible finite conformal module over
$(\ona{Vect}\C^{\ti}, \{L\})$ is isomorphic to one of the modules
$F(\De, \alpha)$ with $\De\neq 0.$

{\em (b)}
Any non-trivial irreducible finite conformal module over the current
algebra
\linebreak
$\gk[t, t^{-1}]$, where $\gk$ is a finite-dimensional
semi-simple Lie algebra, is isomorphic to one of the modules $U[t,
t^{-1}],$ where $U$ is a finite-dimensional non-trivial irreducible
$\gk$-module.
}

\bno
\textbf{Remark 5.2.}\quad Complete reducibility of finite modules over
simple finite Lie conformal algebras breaks down. A classification of
all extensions between finite irreducible modules over semisimple Lie
conformal algebras is given in \cite{CKW}. In particular, it is shown
there that there exists a non-trivial extension of modules
$0\rightarrow M(\De', \alpha')\rightarrow E \rightarrow M(\De, \alpha)\rightarrow 0$
iff $\alpha=\alpha'$ and the pair $(\De, \De')$ takes one of the following
values (cf. \cite{Fuchs}):

\noindent
(i)\quad $\;(a, a),\; (a+2, a),\;(a+3, a),\; (a+4, a)$ where $a\in \C,$

\noindent
(ii)\quad $(1, 0),\;(5, 0),\;(1, -4),\;((7\pm\sqrt{19})/2, (-5\pm
\sqrt{19})/2).$

\bno
\textbf{Remark 5.3.}\quad Theorem 5.1(b) still holds if $\gk$ is a
finite-dimensional simple Lie superalgebra whose maximal reductive
subalgebra is semisimple. However, in the remaining cases, namely
$A(m, n)$ with $m\neq n,$ $C(n)$ and $W(n),$ the description of finite
irreducible $\ona{Cur}\gk$-modules is much more interesting (see
\cite{CK}).

\medskip
The following theorem is a conformal analogue of the classical Lie and
Engel theorems.

\bno
\textbf{Theorem 5.2} \cite{DK}\textbf{.}\quad {\it
{\em (a)}
For any finite module $M$ over a finite solvable Lie conformal algebra
$R$ there exists a non-zero vector $v\in M$ such that
$ a_{\lambda}v=c(a, \lambda)v $, $ a\in R $, $ c(a, \lambda)\in\C[\lambda]$.

{\em (b)} Let $M$ be a finite module over a finite Lie conformal algebra
such that the operator $a_{\lambda}$ is nilpotent on $\C[\lambda]\ot_{\C}M$
for each $a\in R.$ Then there exists a non-zero vector in $M$
annihilated by all operators $a_{\lambda}.$
}

\medskip
It is not difficult to prove the following analogue of the classical
Burnside theorem.

\bno
\textbf{Theorem 5.3.}\quad
\emph{Any non-trivial irreducible finite module over the associative
conformal algebra $\ona{Cur}(\ona{Mat}_N)$ is isomorphic to
$M(\C^N),$ where $\C^N$ is the standard $\ona{Mat}_N$-module. Any
finite module over $\ona{Cur}(\ona{Mat}_N)$ is a direct sum of
irreducible modules.
}

\section{Conformal algebras $\ona{Cend}M$ and $gcM$}

Let $U$ and $V$ be two $\C[\partial]$-modules.
A \emph{conformal linear map} from $U$ to $V$ is a $\C$-linear map
$
  a:U\rightarrow \C[\lambda]\ot_{\C}V
$,
denoted by $a_{\lambda}:U\rightarrow V,$ such that
$\partial a_{\lambda}-a_{\lambda}\partial=-\lambda a_{\lambda}.$
Denote the vector space of all such maps by $\ona{Chom}(U, V).$ It has
a canonical structure of a $\C[\partial]$-module:
$(\partial a)_{\lambda}=-\lambda a_{\lambda}.$

\bno
\textbf{Remark 6.1.}\quad If $U$ and $V$ are modules over a Lie
conformal algebra $R,$ then the $\C[\partial]$-module
$\ona{Chom}(U, V)$ carries an $R$-module structure as well, defined by
$(a_{\lambda}\va)_{\mu}u=a_{\lambda}(\va_{\mu-\lambda}u)
-\va_{\mu-\lambda}(a_{\lambda}u),$ where
$a\in R, u\in U,\;
\va\in \ona{Chom}(U, V).$
Hence we may define the contragredient $R$-module: $U^*= \ona{Chom}(U,
\C),$ where $\C$ is the trivial $R$-module, and the tensor product of
$R$-modules: $U\ot V = \ona{Chom}(U^*, V).$

\medskip
In the special case $U=V=M$ we let $\ona{Cend}M=\ona{Chom}(M, M).$
This $\C[\partial]$-module has the following $\lambda$-product making it an
associative conformal algebra $(a, b\in\ona{Cend}M):$
\setcounter{equation}{22}
\begin{equation}
(a_{\lambda}b)_{\mu}v=a_{\lambda}(b_{\mu-\lambda}v),\quad v\in M.
\end{equation}
Indeed, (12) is immediate, while (13) follows from (23) by replacing
$\mu$ by $\mu+\lambda.$

\bno
\textbf{Remark 6.2.}\quad The $\lambda$-bracket
$[a_{\lambda}b]=a_{\lambda}b-b_{-\lambda-\partial}a$ makes $\ona{Cend}M$ a Lie
conformal algebra denoted by {gcM}. (A decomposition of $M$ in a direct sum of
$\C[\partial]$-modules $M=M_{\bar{0}}\op M_{\bar{1}}$ induces a
$\Z_2$-gradation on the algebra $\ona{Cend}M$ and the $\lambda$-bracket
(16) makes it a Lie conformal superalgebra.) Due to Remark 5.1 one has a
simpler form of this $\lambda$-bracket:
\begin{equation}
[a_{\lambda}b]_{\mu}v=[a_{\lambda}, b_{\mu-\lambda}]v,\quad v\in M.
\end{equation}

\bno
\textbf{Remark 6.3.}\quad Formulas (23) (resp. (24)) show that a
structure of a module over an associative (resp. Lie) conformal
algebra $R$ is the same as a homomorphism of $R$ to the conformal
algebra $\ona{Cend}M$ (resp. $gcM$).

\medskip
Let $N$ be a positive integer, and let
$\ona{Cend}_N\!=\!\ona{Cend}\C[\partial]^N,$ $gc_N\!=\!gc\,\C[\partial]^N,$ where
$\C[\partial]^N$ denotes the free $\C[\partial]$-module of rank $N.$ Remark 6.2
shows that the conformal algebras $\ona{Cend}_N$ and $gc_N$ play the
same role in the theory of conformal algebras as $\ona{End}_N$ and
$gl_N$ play in the theory of associative and Lie algebras. Below we
give a less abstract description of these conformal algebras.

\medskip
Let Diff$_N\C^{\ti}$ be the associative algebra of all $N\ti N$-matrix
valued regular differential operators on $\C^{\ti}.$ It is spanned
over $\C$ by differential operators $At^j\partial^m_t,$ where
$A\sub\ona{Mat}_N,$ $j\in\Z,$ $m\in\Z_+.$ Introduce the following
formal distribution for $A\in\ona{Mat}_N,$ $m\in\Z_+:$
$J_A^m(z)=\sum_{j\in\Z}At^j(-\partial_t)^mz^{-j-1}.$
Then we have:
\begin{equation}
  J_A^m(z)J_B^n(w)
  =
  \sum_{j=0}^m
    \sum_{i=0}^j
      \binom{m}{j}
      \binom{j}{i}
      \partial_w^{j-i}
      J_{AB}^{m+n-j}(w)\partial_w^i\delta(z-w).
\end{equation}

It follows that $F=\{J_A^m(z)\}_{m\in\Z_+,\,A\in\ona{Mat}_N}$ is a
local family, hence
$(\ona{Diff}_N\C^{\ti}, F)$ is a formal
distribution associative algebra. The corresponding associative
conformal algebra is
$$
  \ona{Con}(\ona{Diff}_N\C^{\ti},\, F)
  =
  \sum_{\stackrel{\scriptstyle m\in\Z_+}{A\in\ona{Mat}_N}}\C[\partial]J_A^m
$$
with the $\lambda$-product, derived from (25) to be:
\begin{equation}
\label{eq:26}
  J_{A\;\;\lambda}^mJ_B^n
  =
  \sum_{j=0}^m
    \binom{m}{j}
    (\lambda+\partial)^j
    J_{AB}^{m+n-j}.
\end{equation}

The obvious representation of Diff$_N\C^{\ti}$ on the space $\C^N[t,
t^{-1}]e^{- \alpha t}$ is an irreducible conformal module with the family
$
  E
  =
  \{
    a(z)
    =
    \sum_{m\in\Z}(at^m e^{- \alpha t})z^{-m-1}
  \}_{a\in\C^N}
$.
The associated
module over Con(Diff$_N\C^{\ti}, F)$ is
$\C[\partial]^N=\C[\partial]\ot_{\C}\C^N$ with the $\lambda$-action
\begin{equation}
\label{eq:27}
  J_{A\;\;\lambda}^mv
  =
  (\partial+\lambda + \alpha)^m Av,\quad m\in\Z_+,\;v\in\C^N.
\end{equation}
By Remark 6.3, representation (27) gives us associative conformal
algebra homomorphisms
$\va_{\alpha}:\ona{Con}(\ona{Diff}_N\C^{\ti}, F)\rightarrow \ona{Cend}_N$.

\bno
\textbf{Proposition 6} \cite{DK}\textbf{.}\quad {\it
The homomorphisms $\va_{\alpha}$ are isomorphisms.
}

\medskip\noindent
\emph{Proof}.\quad We have by (27): $(\partial^kJ_A^m)_{\lambda}v
=(-\lambda)^k(\lambda+\partial + \alpha)^m Av.$ Since a conformal linear map is determined
by its values on a set of generators of a $\C[\partial]$-module and the
polynomials $\lambda^k(\lambda+\partial + \alpha)^mv$ with $k, m\in\Z_+,$ $v\in\C^N,$ span
$\C^N[\lambda, \partial],$ the proposition follows.

\medskip
The representation (27) of the associative conformal algebra
$ \ona{Cend}_N $
in $\C[\partial]^N$ can be generalized by keeping formula~(27),
but replacing $\C [\partial]^N$ by $\C [\partial]^N \otimes U$
and $\alpha \in \C$ by an indecomposable linear operator $\alpha$ 
on $U$.  Denote this representation of $\ona{Cend}_N$ by
$\sigma^{as}_{\alpha}$.  This representation gives rise to the representation of the associated Lie
conformal algebra $gc_N$ (see Remark 6.2), which we denote by $\sigma_{\alpha}$.
The representation $\sigma^*_{\alpha}$ contragredient to $\sigma_{\alpha}$ (see Remark 6.1)
is again a free $\C[\partial]$-module with the following $\lambda$-action:
$
  J_{A\;\;\lambda}^mv=-(-\partial- \alpha^*)^m(^tAv)
$,
for
$ m\in\Z_+ $,
$ v\in\C^N \otimes U^*$.

\bno
\textbf{Theorem 6.1} \cite{BKRW}\textbf{.}
\emph{The representations $\sigma_{\alpha}$ and $\sigma^*_{\alpha}$ are all 
non-trivial finite indecomposable $gc_N$-modules.}

\vspace{1ex}
The proof of this theorem relies on methods developed in \cite{KR} and
\cite{CKW}.
The analogue of Theorem 6.1(a),(b) in the associative case is the
following Burnside type theorem.

\bno
\textbf{Theorem 6.2.}\quad
{\it
The $\ona{Cend}_N$-modules $\sigma^{as}_{\alpha}$ are all non-trivial
finite irreducible $\ona{Cend}_N$-modules.}

\vspace{1ex}

Given a collection $P=(P_1(\partial_t), \ldots ,
P_n(\partial_t))$ of non-zero polynomials, one has a formal
distribution subalgebra of $ \mathrm{Diff}_N\C^{\ti}$ consisting
of matrices whose $i$-th column is divisible by $P_i$, $i=1,
\ldots, n$.  The corresponding subalgebra of
$\ona{Cend}_{N}$, denoted by $\ona{Cend}_{N, P}$, still acts
irreducibly on $\C[\partial]^N.$

\bno
\textbf{Conjecture 6.}\quad \emph{If $R\sub\ona{Cend}_N$ is a subalgebra
which still acts irreducibly on $\C[\partial]^N$, then either $R$ 
is conjugate to $\Cur (\Mat_N \C)$ or $R$ is one of
the subalgebras $\ona{Cend}_{N, P}.$
}

\section{$\Ga$-twisted and $\Ga$-local formal distribution algebras
$[$GK$],[$BDK$]$}

We discuss here two generalizations of the notion of a formal
distribution algebra which incorporate the examples of twisted
current algebras and Ramond type superalgebras.

The first generalization requires consideration of non-integral powers
of indeterminates. Let $\Ga$ be an additive subgroup of $\C$
containing $\Z.$ For $\alpha\in\Ga$ we denote by $\bar{\alpha}$ the coset
$\alpha+\Z$.

An $(\bar{\alpha}, \bar{\beta},\ldots)$-\emph{twisted} $U$-valued
formal distribution is a series of the form
$$
  a(z, w, \ldots)
  =
  \sum_{\stackrel{\scriptstyle m\in\bar{\alpha}}{n\in\bar{\beta}}}a_{m, n, \ldots}
  z^{-m-1}w^{-n-1}\ldots.
$$
An $(\bar{\alpha}, \bar{\beta})$-twisted formal distribution $a(z, w)$
is called \emph{local} if, as before,
\hfill \\
$(z-w)^Na(z, w)=0$ for a
sufficiently large $N$.  An example of a $(\bar{\alpha},
-\bar{\alpha})$-twisted local $\C$-valued formal distribution is the
$\bar{\alpha}$-twisted formal $\delta$-function:
\hfill \\
$\delta_{\bar{\alpha}}(z-w)=z^{-1}\sum_{n\in\bar{\alpha}}(w/z)^n.$ As
before, a $(\bar{\alpha}, \bar{\beta})$-twisted local formal
distribution $a(z, w)$ uniquely decomposes in a finite sum of the form
\begin{equation}
a(z, w)=\sum_{j\in\Z_+}c^j(w)\partial^{(j)}_w\delta_{\bar{\alpha}}(z-w).
\end{equation}
Here $c^j(w)$ are $\bar{\alpha}+\bar{\beta}$-twisted formal distributions.

Let now $U$ be an algebra. The pair consisting of an
$\bar{\alpha}$-twisted formal distribution
$a(z)$ and a $\bar{\beta}$-twisted distribution $b(z)$ is called
a \emph{local pair} if
the $(\bar{\alpha}, \bar{\beta})$-twisted formal distribution $a(z)b(w)$ is
local. The $j$-th coefficient of the expansion (28) of $a (z)b(w)$ is
denoted by $a(w)_{(j)}b(w)$ and is called the $j$-th product of these
formal distributions. As before, we define the $\lambda$-product by
$a(w)_{\lambda}b(w)=\sum_{j\in\Z_+}\lambda^{(j)}(a(w)_{(j)}b(w)).$

A pair $(U, F)$, where $U$ is an algebra and $F$ is a family of
twisted $U$-valued formal distributions in $z,$ is called a $\Ga$-{\it
twisted formal distribution algebra} if $\bar{F}$ is a closed
under $j$-th products local family
of $\alpha$-twisted formal distributions with $\alpha\in\Ga,$ whose
coefficients span $U.$ Note that an analogue of Lemma 2 holds.
One can show that Proposition 2 holds if $(U, F)$ is a $\Ga$-twisted
formal distribution algebra. Then $R=\ona{Con}(U, F):=\C[\partial_z]
\bar{F}$ is a conformal algebra (which is associative, commutative or
Lie if $U$ is associative, commutative or Lie respectively). Note also
that $R$ carries a $\Ga/\Z$-gradation
\begin{equation}
R=\op_{\bar{\alpha}\in\Ga/\Z}R_{\bar{\alpha}},
\end{equation}
where $R_{\bar{\alpha}}$
consists of all $\bar{\alpha}$-twisted formal distributions from
$R.$ (It is clear that $R_{\bar{\alpha}(j)}R_{\bar{\beta}}\sub
R_{\bar{\alpha}+\bar{\beta}}$ and $\partial R_{\bar{\alpha}}\sub R_{\bar{\alpha}}.$)

Thus, we get a functor Con from the category of $\Ga$-twisted formal
distribution algebras to the category of $\Ga/\Z$-graded conformal
algebras. Conversely, given a $\Ga/\Z$-graded conformal algebra (29),
we construct the corresponding $\Ga$-twisted formal distribution
algebra by letting (cf. Section 3):
$$
  \ona{Alg}R=\op_{\bar{\alpha}\in\Ga/\Z}(\ona{Alg}R)_{\bar{\alpha}},
$$
where $(\ona{Alg}R)_{\bar{\alpha}}$ is a quotient of the vector space
with the basis $\{a_n\}_{a\in R_{\bar{\alpha}},\,n\in\bar{\alpha}}$ by the
linear span of
$$
  \left\{
    (a+b)_n - a_n - b_n, \,
    (\lambda a)_n - \lambda a_n, \,
    (\partial a)_n + n a_{n-1}
  \right\}_{
    a, b\in R_{\bar{\alpha}}, \; n\in\bar{\alpha}
  }
$$
with the product:
$$
  a_mb_n
  =
  \sum_{j\in\Z_+}
    \binom{m}{j}
    (a_{(j)}b)_{m+n-j},
  \qquad
  \textrm{where} \; \;
  a\in R_{\alpha}, \;
  b\in R_{\beta}, \;
  m\in\bar{\alpha}, \;
  n\in\bar{\beta}.
$$
The local family is $F=\cup_{\bar{\alpha}\in\Ga/\Z}
\{a(z)=\sum_{n\in\bar{\alpha}}a_nz^{-n-1}\}_{a\in R_{\bar{\alpha}}}.$
(Of course in the non-twisted case ($\Ga=\{1\}$) we recover the
construction of Section 3.) The relations between the functors Con and
Alg are the same as in the untwisted case (see Section 3).

\bno
\textbf{Example 7.1.}\quad Let $A=\op_{\bar{\alpha}\in\Ga/\Z}
A_{\bar{\alpha}}$ be a $\Ga/\Z$-graded algebra. Consider the following
subalgebra of the algebra $\op_{\alpha\in\Ga}At^{\alpha}:$
$U=\sum_{\alpha\in\Ga}A_{\bar{\alpha}}t^{\alpha}.$
This is a $\Ga$-twisted formal distribution algebra with the local
family $F=\cup_{\bar{\alpha}\in\Ga/\Z}F_{\bar{\alpha}},$ where
$F_{\bar{\alpha}}=\{a(z)=\sum_{n\in\bar{\alpha}} (at^n)z^{-n-1}\}_{a\in
A_{\bar{\alpha}}}.$ (Indeed, we have:
$a(z)b(w)=(ab)(w)\delta_{\bar{\alpha}}(z-w)$
for $a\in A_{\bar{\alpha}},$ $b\in A_{\bar{\beta}}.$)
It is called a $\Ga$-\emph{twisted current algebra}. Thus, we have:
$\ona{Con}(U, F)=\op_{\bar{\alpha}\in\Ga/\Z}\C[\partial]A_{\bar{\alpha}}$ is a
$\Ga/\Z$-graded current conformal algebra.

\bno
\textbf{Remark 7.1.}\quad The classification of gradings of $R$ depends
on the description of its automorphism group Aut$R$. It is easy to see
that
$\ona{Aut}(\ona{Cur}A)=\ona{Aut}A,$
provided that $aA\neq 0$ if $a\neq 0.$ Indeed, applying an
automorphism $\sigma$ to $a_{\lambda}b=ab,$ $a, b\in A,$ we have $P(-\lambda)
Q(\lambda+\partial)ab= R(\partial)ab,$ where $\sigma(a)=P(\partial)a,$ etc. It follows that
$P(\lambda)$ is independent of $\lambda,$
hence $\sigma(A)\sub A.$ One can show in a similar way
that the group of automorphisms of the Virasoro conformal algebra is
trivial.

\bno
\textbf{Remark 7.2.}\quad Let $\sigma$ be an order $m$ automorphism of a
conformal algebra $R$. It defines a $\Ga=\frac{1}{m}\Z/\Z$-grading of
$R$. Let Alg$(R, \sigma)$ denote the corresponding maximal $\Ga$-twisted
formal distribution algebra. One can show that Alg$(R, \sigma_i),$ $i=1,
2,$ are isomorphic if $\sigma_1$ and $\sigma_2$ lie in the same connected
component of the group Aut$R$.

\medskip
The second generalization deals with the usual formal distributions,
but a more general notion of locality. Let $\Ga$ be a multiplicative
subgroup of $\C^{\ti}.$ A formal distribution $a(z, w)$ is called
$\Ga$-local if
$P(z/w)a(z, w)=0$
for some polynomial $P(x)$ all of whose roots lie in $\Ga.$ Obviously,
for $\Ga=\{1\}$ we have the usual locality. The special case when
$P(x)$ has no multiple roots was studied in detail in \cite{GK}. An
example of a $\Ga$-local $\C$-valued formal distribution is $\delta(z-\alpha
w)$ where $\alpha\in\Ga.$ As before, a $\Ga$-local $U$-valued formal
distribution $a(z, w)$ uniquely decomposes in a finite sum of the
form:
\begin{equation}
  a(z, w)
  =
  \sum_{\stackrel{\scriptstyle j\in\Z_+}{\alpha\in\Ga}}
    c^{j, \alpha} (w)\partial_w^{(j)}\delta(z-\alpha w),
\end{equation}
where $c^{j, \alpha}(w)$ are some $U$-valued formal distributions.

Let now $U$ be an algebra. The pair of
$U$-valued formal distributions $a(z)$
and $b(z)$ is called $\Ga$-\emph{local} if the formal distribution
$a(z)b(w)$ is $\Ga$-local. The $(j, \alpha)$-coefficient of the expansion
(30) of $a(z)b(w)$ is denoted by $a(w)_{(j, \alpha)}b(w)$ and is called
the $(j, \alpha)$-th product of $a(w)$ and $b(w)$. As before, we define
the $(\lambda, \alpha)$-product $(\alpha\in\Ga)$ by
$a(w)_{\lambda, \alpha}b(w)=\sum_{j\in\Z_+}\lambda^{(j)}(a(w)_{(j, \alpha)}b(w))$
and let $a(w)_{\lambda, \Ga}b(w)= a(w)_{\lambda, 1}b(w).$

For $\alpha\in\Ga$ introduce the following linear operator $T_{\alpha}$ on
the space of formal distributions:
$
  T_{\alpha}(a(z))=\alpha a(\alpha z)
$.
Then one has:
\begin{eqnarray}
  a(w)_{\lambda, \alpha}b(w)
  & = &
  (T_{\alpha}a(w))_{\lambda, \Ga}b(w),\\
  T_{\alpha}(a(w)_{\lambda, \Ga}b(w))
  & = &
  (T_{\alpha}a(w))_{\alpha\lambda,
  \Ga}(T_{\alpha}b(w)).
\end{eqnarray}

A pair $(U, F)$, where $U$ is an algebra and $F$ is a family of
$U$-valued formal distributions in $z,$ is called a $\Ga$-\emph{local
formal distribution algebra} if $\bar{F}$ consists of pairwise
$\Ga$-local formal distributions whose coefficients span $U.$ Here
$\bar{F}$ denotes the minimal family containing $F$, invariant under
all $T_{\alpha}\;(\alpha\in\Ga)$ and closed under all $(j, \alpha)$-products
with $\alpha\in\Ga,$ $j\in\Z_+.$ Again, an analogue of Lemma 2 holds.
Again, one can show that Proposition 2 holds for the $(\lambda,
\Ga)$-product. Thus, a $\Ga$-local formal distribution algebra $(U, F)$
gives rise to a conformal algebra
$R=\ona{Con}(U, F)=\C[\partial_z]\bar{F},$ and an action of $\Ga$ on it
by \emph{semilinear} automorphisms.
``Semilinear''
means that we have a homomorphism $\alpha\mapsto T_{\alpha}$ of $\Ga$ to the
group of $\C$-linear invertible maps of $R$ such that:
$\partial T_{\alpha}=\alpha T_{\alpha}\partial,$
$T_{\alpha}(a_{(j)} b)=\alpha^j(T_{\alpha}a)_{(j)}(T_{\alpha}b).$
Furthermore, for each pair $a, b\in R,$ one has:
\begin{equation}
  (T_{\alpha}a)_{(j)}b
  =
  0
  \quad
  \mathrm{for~all~but~finitely~many~} \alpha\in\Ga, j\in\Z_+.
\end{equation}

Conversely, given a conformal algebra $R$ with an action of a group
$\Ga\sub\C^{\ti}$ by semilinear automorphisms, such that (33) holds
we construct the
corresponding $\Ga$-local formal distribution algebra Alg$(R, \Ga)$
which as a vector space is the quotient of $R[t, t^{-1}]$ by
the linear span of elements (as before $a_n$ stands for $at^n$):
$\{(\partial a)_n+na_{n-1},$ $(T_{\alpha}a)_n-\alpha^{-n}a_n\}_{a\in
R,\,n\in\Z},$
with the following product (cf. \cite{GK} and Section 3):
$$
  a_m b_n
  =
  \sum_{\stackrel{\scriptstyle j\in\Z_+}{\alpha\in\Ga}}
    \alpha^m
    \binom{m}{j}
    ((T_{\alpha}a)_{(j)}b)_{m+n-j}.
$$
The $\Ga$-local family is
$
  F=\{a(z)=\sum_{m\in\Z}a_mz^{-m-1}\}_{a\in R}
$.

\bno
\textbf{Example 7.2.} (cf. \cite{GK})\quad
Let $A$ be an algebra with an action of a finite group
$\Ga\sub\C^{\ti}$ by automorphisms. The action of $\Ga$ on $A$ extends
to Cur$A=\C[\partial]\ot_{\C}A$ using $\partial T_{\alpha}=\alpha T_{\alpha}\partial,$ $\alpha\in
\Ga.$ The corresponding $\Ga$-local formal distribution algebra is:
$a(z)b(w)=\sum_{\alpha\in\Ga}((T_{\alpha}a)b)(w)\delta(z-\alpha w),$
where $a, b\in A.$
It is easy to see that we get once more a $\Ga$-twisted current
algebra.

\bno
\textbf{Remark 7.3.}\quad The simplest case of simple poles considered
in \cite{GK} is the case of $\Ga$-local formal distribution algebras
which correspond to conformal algebras with the trivial action of
$\partial$ and the $\lambda$-product independent of $\lambda,$ i.e. ordinary
algebras (with an action of $\Ga$). The so called sin algebra \cite{GL}, which
is a $q$-analogue of Diff$_N\C^{\ti},$ is a $\Ga$-local formal
distribution algebra associated to the algebra of infinite matrices,
where $\Ga=\{q^n\}_{n\in\Z}.$ The role of this algebra is analogous to
that of $\ona{Cend}_N$ in the theory of ordinary associative conformal
algebras \cite{GK}.

\bno
\textbf{Remark 7.4.}\quad One defines $\Ga$-twisted modules and
$\Ga$-conformal modules in the obvious way, and establishes
equivalence to $\Ga$-graded and $\Ga$-equivariant conformal modules
respectively.

\section{Work in progress}

\noindent
\textbf{8.1. Case of several indeterminates}
\cite{BDK}\textbf{.}\quad Let $z=(z_1, \ldots, z_n),$
\hfill \\
$w=(w_1, \ldots, w_n),$ $\partial=(\partial_1, \ldots, \partial_n)$ and
$\lambda=(\lambda_1, \ldots, \lambda_n).$ A formal distribution
$a(z, w)$ with values in an algebra $U$ is called {\it local} if
$$
  (z_i-w_i)^Na(z, w)=0,
  \mathrm{~for~} i=1, \ldots, n
  \mathrm{~and~} N \gg 0.
$$
We have a finite expansion similar to (1):
$$
  a(z, w)
  =
  \sum_{j\in\Z_+^n}c^j(w)\partial^{(j)}_w\delta(z-w),
$$
where
$$
  \partial^{(j)}_w\delta(z-w)=\prod_i\partial_{w_i}^{(j_i)}\delta(z_i-w_i).
$$
All definitions and statements of Sections 1-3 extend without difficulty
to the case of several indeterminates.

\bno
\textbf{Example 8.1.}

\noindent
(a) The Lie algebra $W_n$ of all derivations
of the algebra $\C[x_1, x_1^{-1}, \ldots, x_n, x_n^{-1}]$ is spanned
by pairwise local formal distributions
$
  A^i(z)
  =
  -\delta(z-x)\partial/\partial x_i
$.
The associated conformal algebra (in $n$ indeterminates) is
$\ona{Con}W_n=\sum_{i=1}^n\C[\partial]A^i,$
with the $\lambda$-bracket
$$
  [A^i_{\lambda}A^j]=\partial_iA^j+\lambda_iA^j+\lambda_jA^i.
$$

\noindent
(b) The subalgebra $S_n$ of divergence 0 derivations is a formal
distribution subalgebra of $W_n.$ The corresponding conformal algebra
is
\hfill \\
$\ona{Con}S_n=$
$\left\{\sum_iP_i(\partial)A^i\left|\sum_iP_i\partial_i=0
\right.\right\}.$

\noindent
(c) The subalgebra $H_n,$ $n=2k,$ of Hamiltonian derivations is a
formal distribution subalgebra of $W_n.$ The corresponding conformal
algebra is:
$\ona{Con}H_n=\C[\partial]A$
with the $\lambda$-bracket
$$
  [A_{\lambda}A]=\sum_{i=1}^k(\lambda_{k+i}\partial_iA-\lambda_i\partial_{k+i}A).
$$

\noindent
(d) The subalgebra $K_n,$ $n=2k+1,$ of contact derivations is not a
finite over $\C[\partial_1, \ldots,\partial_n]$ formal distribution
algebra in $n$ indeterminates.

\medskip
The structure theory of conformal algebras in several indeterminates
is being worked out.

\bno
\textbf{Remark 8.1.}\quad The algebra $K_n$ is finite (with one
generator) over the (non-commu\-tative) Weyl algebra in $k$
indeterminates. This leads to a more general
notion of a conformal algebra where
$\C[\partial_1, \ldots,\partial_n]$ is replaced by a bialgebra. This is
yet a special case of a general notion of a Lie$^*$ algebra introduced
in \cite{BD}. The structure theory of finite Lie$^*$ algebras is being worked out in \cite{BDK}.

\bno
\textbf{8.2. Cohomology} \cite{BKV}\textbf{.}\quad
Let $R$ be a Lie conformal algebra and let $M$ be an $R$-module. An
$n$-\emph{cochain of $R$ with coefficients in} $M$ is a skewsymmetric
conformal anti-linear map (cf. Section 6), i.e. a $\C$-linear map
$$
  \ga:R^{\ot n} \longrightarrow \C[\lambda_1, \ldots, \lambda_n]\ot M,\quad
  a_1\ot\cdots\ot a_n \mapsto\ga_{\lambda_1, \ldots, \lambda_n}(a_1, \ldots,
  a_n)
$$
such that
$
  \ga_{\lambda_1, \ldots, \lambda_n}
  (a_1, \ldots, \partial a_i, \ldots, a_n)
  =
  -\lambda_i\ga_{\lambda_1, \ldots, \lambda_n}
  (a_1, \ldots, a_i, \ldots, a_n)
$,
and $\ga$ is
\hfill \\
skewsymmetric with respect to simultaneous permutations of
$a_i$'s and $\lambda_i$'s.

A differential $d\ga$ of a cochain $\ga$ is defined by the following
formula:
\begin{eqnarray*}
  &
  \! \! \! \!
  (d\ga)_{\lambda_1, \ldots, \lambda_{n+1}}(a_1, \ldots, a_{n+1})
  =
  -{\displaystyle \sum_{i=1}^{n+1}}
    (-1)^ia_{i_{\lambda_i}}\ga_{\lambda_1,\ldots,\hat{\lambda}_i,
  \ldots,\lambda_{n+1}}(a_1,\ldots,\hat{a}_i,\ldots,a_{n+1})\\
  &
  +
  {\displaystyle \sum_{\stackrel{\scriptstyle i, j = 1}{i < j}}^{n+1}}
    (-1)^{i+j}
    \ga_{\lambda_i+\lambda_j, \lambda_1, \ldots, \hat{\lambda}_i,
      \ldots, \hat{\lambda}_j, \ldots, \lambda_{n+1}}
  ([a_{i_{\lambda_i}}a_j], a_1, \ldots, , \hat{a}_i, \ldots, \hat{a}_j,
  \ldots, a_{n+1}).
\end{eqnarray*}
Then $d\ga$ is again a cochain and $d^2\ga=0.$ Thus the cochains form
a complex, called the \emph{basic complex} and denoted by
$\widetilde{C}^*(R,M)=\op_{n\in\Z_+}\widetilde{C}^n(R, M).$

Define a structure of a $\C [\partial]$-module on $\widetilde{C}^*
(R,M)$ by
$$
  (\partial \ga)_{\lambda_1, \ldots , \lambda_n}
  =
  (\partial + \lambda_1 + \ldots + \lambda_n )
  \ga_{\lambda_1 , \ldots , \lambda_n} \, .
$$

Then $\partial$ commutes with $d$ and therefore we can define the
\emph{reduced} complex
$$
  C^*(R, M)
  =
  \widetilde{C}^*(R, M)/ \partial \widetilde{C}^*(R,M) \, .
$$
We define the \emph{cohomology}
$H^*(R, M)=\op_{n\in\Z_+} H^n(R, M)$ of $R$ with coefficients in a
module $M$ (resp. $\widetilde{H}^* (R,M)$) to be the cohomology
of the reduced complex (resp. of the basic complex).

As in the case of Lie algebra cohomology, $\widetilde{H}^0(R, M)={M}^R$,
\hfill \\
$H^1(R, Chom(N,M))$ parameterizes extensions of $M$ by a module $N$,
$H^2(R, \C)$ parameterizes central extensions of $R$, $H^2(R, M)$
parameterizes abelian extensions of $R$, etc.

\bno
\textbf{Example 8.2.}

\noindent
(a)\quad Let $R$ be the Virasoro conformal algebra and $M=\C$ be its
trivial module. Then $\widetilde{H}^n (R,\C)$ is
1-dimensional for $n=0$ or 3 and is 0 otherwise. It follows that
$H^n(R, \C)$ is 1-dimensional for $n=0, 2$ or 3 and is 0 otherwise.
This example is intimately related to the Gelfand-Fuchs cohomology (see
\cite{Fuchs}).  We also have calculated $H^* (R,M (\De , \alpha ))$.

\noindent
(b)\quad Let $R=\ona{Cur}\gk$, where $\gk$ is a simple
finite-dimensional Lie algebra. Then $\widetilde{H}^* (R, \C)$ is the Grassmann algebra $G=\op_jG_j$ on generators of degrees
$2m_1+1, 2m_2+1, \ldots,$ where $m_i$'s are the exponents of $\gk$, i.e. it is
the same as $H^*(\gk, \C).$ It follows that $H^*(R, \C)=(\op_jG_j)\op
(\op_jG_{j-1}).$ In particular $\dim H^2(R, \C)=1.$  We also have
calculated $H^* (R,M (U))$.

\medskip
In the case when $R$ is an associative conformal algebra one can
construct analogues of Hochschild and cyclic cohomology.

\bno
%\textbf{8.3. Infinite irreducible subalgebras of $gc_N$.}\quad
%Let $f(u, v)$ be a non-degenerate $\C$-valued bilinear form on
%$\C[\partial]^N$ such that $\partial^*=-\partial.$
%Let $oc_N$ (resp. $spc_N$)$=\{a\in gc_N|f(a_{\lambda}u, v)+
%f(u, a_{\lambda}v)=0,$ $u, v\in\C[\partial]^N\}$ if $f$ is symmetric
%(resp. skewsymmetric). These are \emph{orthogonal} (resp. {\it
%symplectic}) Lie conformal algebras. They are infinite conformal
%subalgebras of $gc_N$ which still act irreducibly on $\C[\partial]^N.$

%\bno
\noindent\textbf{8.3.~Open Problem. }\quad{\it
Classify all infinite conformal subalgebra of $gc_N$ which still acts
irreducibly on $\C[\partial]^N$.
}

\bno
\textbf{Remark 8.3.}\quad All finite subalgebras of $gc_N$ that act
irreducibly on $\C[\partial]^N$ are described in \cite{DK}.

\section*{Acknowledgements}

Partially supported by NSF grant DMS-9622870.

%\section*{References}

\end{document}